\documentclass[12pt]{article}
\usepackage{amsmath,amsthm,amsfonts, latexsym, amssymb}
\usepackage{amscd}
\usepackage{enumerate}

\def\CC{\mathbb C}
\def\DD{\mathbb D}
\def\RR{\mathbb R}

\def\Hol{\mathop{\rm Hol}\nolimits}
\def\re{\mathop{\rm Re}\nolimits}

\def\eps{\varepsilon}
\def\ol{\overline}
\def\pa{\partial}

\newtheorem{thm}{Theorem}[section]

\newtheorem{exam}{Example}[section]
\newtheorem{prop}{Proposition}[section]
\newtheorem{cor}{Corollary}[section]
\newtheorem{rem}{Remark}[section]

\newtheorem{defn}{Definition}[section]

\def\beginpf{\begin{proof}}
\def\endpf{\end{proof}}

\def\beq{\begin{equation}}
\def\eeq{\end{equation}}

\def\esssup{\mathop{\rm ess~sup}\nolimits}
\def\re{\mathop{\rm Re}\nolimits}
\def\im{\mathop{\rm Im}\nolimits}

\begin{document}

\author{\rm Isabelle Chalendar\\ \small{Universit\'e de Lyon; CNRS; Universit\'e Lyon 1; INSA de Lyon; Ecole Centrale de Lyon}\\ \small{ CNRS, UMR 5208, Institut Camille Jordan}\\ \small{ 43 bld. du 11 novembre 1918, F-69622 Villeurbanne Cedex, France}\\ \small{\tt E-mail: chalendar@math.univ-lyon1.fr}
\and Jonathan R. Partington\\\small{School of Mathematics, University of Leeds, Leeds LS2 9JT, U.K.} \\ \small{\tt E-mail: J.R.Partington@leeds.ac.uk}}

\title{Phragm\'en--Lindel\"of principles for generalized analytic functions on unbounded domains}

\maketitle

\begin{abstract}
We prove versions of the Phragm\'en--Lindel\"of strong maximum principle for generalized analytic functions defined on
unbounded domains. A version of Hadamard's three-lines theorem is also derived.
\end{abstract}



\section{Introduction}

Versions of the maximum principle for complex-valued functions defined on a domain in $\CC$ have been
of interest since the development of the classical maximum modulus theorem and
Phragm\'en--Lindel\"of principle for holomorphic functions (see, e.g. \cite[Chap. V]{titchmarsh}).
It is important to distinguish between two types of result here.
First, there is the  {\em weak maximum principle\/}
asserting that under certain circumstances a nonconstant function $f: \Omega \to \CC$  cannot attain a local maximum in its domain
$\Omega$: thus if $\Omega$ is bounded and $f$ is continuous on $\overline\Omega$ we have
\begin{equation}\label{eq:maxom}
\sup_{z \in \Omega} |f(z)| = \sup_{z \in  \partial\Omega} |f(z)|.
\end{equation}

Second -- and this will be our main concern in this paper -- there is the
{\em strong maximum principle\/} or {\em  Phragm\'en--Lindel\"of principle}. This generally applies to unbounded
domains, and generally a supplementary hypothesis on $f$ is required for the conclusion (\ref{eq:maxom}) to hold.
For example, if $f: \Omega \to \CC$ is analytic, where $\Omega=\CC_+$, the right-hand half-plane $\{z \in \CC: \re z > 0\}$,
then if $f$ is known to be bounded we may conclude that (\ref{eq:maxom}) holds, whereas the example
$f(z)=\exp(z)$ shows that it does not hold in general.
\\

We shall use the following standard notation:\\ 
\[
\partial f= \dfrac{\partial f}{\partial z}=\frac12 (f_x-if_y) \quad \hbox{ and} \quad \overline\partial f=\dfrac{\partial f}{\partial \overline z}
=\frac12 (f_x+if_y).
\]

For quasi-conformal mappings $f$, that is, those satisfying the Beltrami equation $\overline\partial f = \nu \partial f$ with
$|\nu| \le \kappa < 1$, the weak maximum principle holds (see, for example \cite{chen}). This fact was used in \cite[Prop. 4.3.1]{BLRR}
to deduce a weak maximum principle for functions solving the conjugate Beltrami equation 
\beq\label{eq:conjbelt}
\overline\partial f = \nu \overline{\partial f}.
\eeq
Their argument is based on the fact that if $f$ is a solution to (\ref{eq:conjbelt}), then
it also satisfies a classical Beltrami equation
$\overline\partial f = \nu_f \partial f$,
where $\nu_f(z)=\nu(z) \overline{\partial f(z)}/\partial f(z)$, and hence $f=G \circ h$ where $G$ is holomorphic and $h$ is a quasi-conformal
mapping (cf. \cite[Thm. 11.1.2]{IM01}).

Carl \cite{carl} considered functions $w$ satisfying equations of the form
\beq\label{eq:carl}
\overline\partial w(z)+A(z)w(z)+B(z)\overline{w(z)} = 0
\eeq
and deduced a weak maximum principle for such functions,
analogous to (\ref{eq:maxom}),  under certain hypotheses on the functions $A$ and $B$.
We shall take this as our starting point.\\

For general background on generalized analytic functions (pseudo-analytic functions) we refer to
the books \cite{bers, kravchenko,vekua}.
The following definitions are taken from the recent paper \cite{BLRR}.

\begin{defn}
Let
$1 \le p < \infty$.
For $\nu \in W^{1,\infty}(\DD)$ (i.e., a Lipschitz function with bounded
partial derivatives), the class $H^p_\nu$ consists of
all measurable functions $f: \DD \to \CC$ satisfying the conjugate Beltrami equation
(\ref{eq:conjbelt})
in a distributional sense, such that the norm
\[
\|f\|_{H^p_\nu}=\left( \esssup_{0<r<1} \frac{1}{2\pi} \int_0^{2\pi} |f(re^{it})|^p \, dt \right)^{1/p}
\]
is finite. Clearly for $\nu=0$ we obtain the classical Hardy space $H^p(\DD)$.
If instead $\nu$ is defined on an arbitrary subdomain $\Omega \subset \CC$,
we may define the class $H^\infty_\nu(\Omega)$ as the space of all bounded measurable functions satisfying  
(\ref{eq:conjbelt}), equipped with the supremum norm.\\

We may analogously define spaces $G^p_\alpha(\DD)$, where $\alpha \in L^\infty(\DD)$, and in general $G^\infty_\alpha(\Omega)$,
where now, for a function $w$ we replace (\ref{eq:conjbelt}) by
\beq\label{eq:dbar}
\overline\partial w = \alpha \overline w.
\eeq
Once again, the case $\alpha=0$ is classical.
\end{defn}

When $\nu$ is real (the most commonly-encountered situation), 
there is a link between the two notions: suppose that $\|\nu\|_{L^\infty(\Omega)}$ with $\|\nu\|_\infty \le \kappa < 1$, and
set $\sigma=\dfrac{1-\nu}{1+\nu}$ and $\alpha= \frac{\overline\partial \sigma}{2\sigma}$, so that $\sigma \in W^{1,\infty}_\RR(\Omega)$.
Then $f \in L^p(\DD)$ satisfies  (\ref{eq:conjbelt}) if and only if
$w:= \dfrac{f-\nu \overline f}{\sqrt{1-\nu^2}}$ satisfies (\ref{eq:dbar}).\\

We shall mainly be considering the class $G^\infty_\alpha$, for which it is possible to prove a strong maximum
principle and a generalization of the Hadamard three-lines theorem under mild hypotheses on $\alpha$, which are
satisfied in standard examples. The referee has suggested that there may be a link between these assumptions and the
strict ellipticity of $\sigma$, although we have not been able to show this.

\section{Functions defined on unbounded domains}

The following result is an immediate consequence of \cite[Thm. 1]{carl}, taking
$A=0$ and $B(z)=-\alpha(z)$ in (\ref{eq:carl}) in order to obtain  (\ref{eq:dbar}).

\begin{prop}\label{prop:carl}
Suppose that $\Omega$ is a bounded domain in $\CC$ and that $w$ is a continuous function on $\overline \Omega$
such that (\ref{eq:dbar}) holds in $\Omega$, where $\alpha$ satisfies $2|\alpha|^2 \ge | \partial \alpha|$.
Then 
$|w(z)| \le \sup_{\zeta \in \partial \Omega}|w(\zeta)|$ for all $z \in \Omega$.
\end{prop}

\beginpf
Taking $k=2$ in \cite[Thm. 1]{carl}, we require that the matrix $M=(m_{ij})_{i,j=1}^2$ be negative semi-definite,
where, with $a=-2|\alpha|^2 $ and $b=-\partial \alpha$, we have
\[
M=\begin{pmatrix}a + \re b & \im b \\ \im b & a-\re b\end{pmatrix}.
\]

On calculating $m_{11}$, $m_{22}$ (which must be non-positive) and $\det M$  (which must be non-negative) we obtain the 
sufficient conditions
$-2|\alpha|^2 \pm \re \partial \alpha \le 0$ and $2|\alpha|^2 \ge | \partial \alpha|$: clearly the second condition
implies the first.
\endpf

\begin{exam} {\rm 
In the example $\sigma=1/x$, occurring in the study of the tokamak reactor
\cite{fl,flps}, we have $\alpha(x)=-\frac{1}{4x}$ and $\partial\alpha=\frac{1}{8x^2}$; thus the inequality
$2|\alpha|^2 \ge | \partial \alpha|$ is always an equality.

Note that by rescaling $z$ we may transform the equation (\ref{eq:dbar}) to
one with $\alpha= -\frac{1}{\lambda x}$ for any $\lambda>0$ (with the domain also changing); then the inequality
requires that $2/\lambda^2 \ge 1/2\lambda$, so that if we take $0<\lambda<4$ the inequality is strict.
}
\end{exam}


Now for $\eps>0$ we write $h_\eps(z)=1/(1+\eps z)$, and note that 
whenever $\Omega \subset \CC_+$ is a domain, we have that
the functions $h_\eps$ satisfy
\begin{enumerate}[(i)]
\item\label{en1} For all $\eps>0$, $h_\eps \in \Hol(\Omega) \cap C(\overline\Omega)$.
\item\label{en2} For all $\eps>0$, $\lim_{|z| \to \infty, z \in \overline\Omega} h_\eps(z)=0$.
\item\label{en3} For all $z \in \Omega$, $\lim_{\eps \to 0} |h_\eps(z)|=1$.
\item\label{en4} For all $\eps>0$, for all $z \in \partial\Omega$, $|h_\eps(z)| \le 1$.
\end{enumerate}

Suppose that $\overline\partial w=\alpha \overline w$ and that $h$ is holomorphic; then
$\overline\partial(hw)=\beta \overline{hw}$, where $\beta=\alpha h/\overline h$. Moreover,
\[
\partial \beta = \partial(\alpha h)/\overline h= (\partial \alpha)(h/\overline h) + \alpha (\partial h)/\overline h.
\]

That is, with $h=h_\eps$, we have $|\beta|=|\alpha|$ and $|\partial \beta| \le |\partial \alpha| +   |\alpha| |\partial h_\eps|/|h_\eps|$.

\begin{thm}\label{thm:halfplane}
Suppose that $\Omega \subset \CC_+$ (not necessarily bounded) and that 
$w$ is a continuous bounded function on $\overline \Omega$
such that (\ref{eq:dbar}) holds in $\Omega$ where $\alpha$ is a $C^1$ function satisfying $2|\alpha|^2 \ge | \partial \alpha|+
|\alpha| |\partial h_\eps|/|h_\eps|$ for all $\eps>0$.
Then 
$|w(z)| \le \sup_{\zeta \in \partial \Omega}|w(\zeta)|$ for all $z \in \Omega$.
\end{thm}

\beginpf
Fix $\eps>0$ and $M=\sup_{\zeta \in \partial \Omega}|w(\zeta)|$. Suppose that $M>0$. Then by property (\ref{en2}) there is an $\eta>0$ such that
for all $z \in \overline\Omega$ with $|z| \ge \eta$ we have $|w(z) h_\eps(z)| \le M$.

Now, by property (\ref{en1}) and Proposition \ref{prop:carl} we have
\[
\sup_{z \in \Omega \cap D(0,\eta)} |w(z)h_\eps(z)| = \sup_{z \in \partial (\Omega \cap D(0,\eta))}|w(z)h_\eps(z)|,
\]
at least if $2|\alpha|^2 \ge | \partial \alpha|+  |\alpha||\partial h_\eps|/|h_\eps| $.\\

Now $\partial (\Omega \cap D(0,\eta)) \subset (\partial\Omega \cap \overline{D(0,\eta)}) \cup (\partial D(0,\eta) \cap \overline\Omega)$.

By hypothesis, $|w(z)| \le M$ if $z \in \partial\Omega$, and by property (\ref{en4}), $|h_\eps(z)| \le 1 $ for $z \in \partial \Omega$.
So $\sup_{z \in \partial\Omega \cap \overline{D(0,\eta)}} |w(z) h_\eps(z)| \le M$.

By the definition of $\eta$ we also have $|w(z) h_\eps(z)| \le M$ if $|z| \ge \eta$ with $z \in \overline \Omega$, and in
particular for $z \in \overline\Omega \cap \partial D(0,\eta)$.

We conclude that $\sup_{z \in \Omega \cap D(0,\eta)} |w(z) h_\eps(z)| \le M$. However, $|w(z) h_\eps(z)| \le M$ whenever 
$z \in \overline\Omega$ with $|z| \ge \eta$, and hence
$\sup_{z \in \Omega} |w(z) h_\eps(z)| \le M$. Now, letting $\eps$ tend to $0$, and using property (\ref{en3}), we have
the result in the case $M>0$.\\

If $M=0$, then by the above we have that $\sup_{z \in \partial \Omega} |w(z)| \le \gamma$ for all $\gamma>0$, and the
same holds for $z \in \Omega$ by the above. Letting $\gamma \to 0$ we conclude that $w$ is identically $0$ on $\Omega$.

\endpf

\begin{exam}\label{ex:2.2}
 {\rm 
Consider the case $\alpha= -\frac{1}{\lambda x}$ and $\partial \alpha = \frac{1}{2\lambda x^2}$. For the
hypotheses of the theorem to be valid we require
\[
\frac{2}{\lambda x^2} \ge \frac{1}{2\lambda x^2}+\frac{1}{\lambda x}\frac{\eps}{|1+\eps z|}.
\]
If $\lambda=1$ (and by rescaling the domain we can assume this) then
this always holds, since   $|1+\lambda z| \ge \lambda x$.
}
\end{exam}


In the following theorem, it will be helpful to note that we shall
be considering composite mappings as follow:
\[
\Lambda \xrightarrow{h} \Omega \xrightarrow{w} \CC \qquad \hbox{and}
\qquad
\Lambda \xrightarrow{h} \Omega \xrightarrow{\alpha} \CC.
\]

\begin{thm}\label{thm:notdense}
Suppose that $\Omega \subset \CC$  is simply-connected and that the disc $D(a,r)$ is contained in $ \CC \setminus \overline{\Omega}$.  
Let $h:  \CC \to \CC$ be defined by $h(z)=re^z+a$, and let
$\Lambda$ be a component of $h^{-1}(\Omega)$.
Set $g_\eps(z)=1/(1+\eps g(z))$, where $g(z)=\log\left (\dfrac{z-a}{r}\right)$ is a single-valued inverse to $h$ defined on $\Omega$.
Suppose that
$w$ is  a continuous bounded function on $\overline \Omega$
such that (\ref{eq:dbar}) holds in $\Omega$ with $\alpha$  a $C^1$ function satisfying 
\beq\label{eq:alpha}
2|\alpha|^2 \ge | \partial \alpha|+
|\alpha| |\partial g_\eps|/|g_\eps|
\eeq
 for all $\eps>0$.
Then 
$|w(z)| \le \sup_{\zeta \in \partial \Omega}|w(\zeta)|$ for all $z \in \Omega$.
\end{thm}


\beginpf
First we identify the equation satisfied by $v=w \circ h$, where $h$ is holomorphic. Namely,
\begin{eqnarray*}
\overline\partial  v &=& \overline\partial (w \circ h )= \overline{\partial(\overline w \circ h)}=\overline{(\partial \overline w \circ h)(\partial h)}=(\overline\partial w \circ h)(\overline{\partial h})\\
&=& ((\alpha\overline w)\circ h)(\overline{\partial h})
= (\alpha \circ h)(\overline w \circ h) (\overline{\partial h})= \beta \overline v,
\end{eqnarray*}
where $\beta=(\alpha \circ h)(\overline{\partial h})$. Note that $\partial\beta = (\partial\alpha \circ h)|\partial h|^2$, since $\partial (\overline{\partial h})=0$.

The condition 
\beq
\label{eq:beta}
2|\beta|^2 \ge |\partial\beta| +|\beta|  |\partial h_\eps|/|h_\eps|
\eeq
 at a point of $\Lambda$ can be rewritten
\[
2|\alpha \circ h|^2 |\partial h|^2 \ge |\partial\alpha \circ h|\,|\partial h|^2 + |\alpha \circ h| \, |\partial h| |\partial h_\eps|/|h_\eps|.
\]

Now  $g_\eps=h_\eps \circ g$; thus $\partial h_\eps=(\partial g_\eps \circ h)(\partial h)$.

That is, (\ref{eq:beta}) is equivalent to
\[
2|\alpha \circ h|^2 |\partial h|^2 \ge |\partial\alpha \circ h|\,|\partial h|^2 + |\alpha \circ h| \, |\partial h|^2 |\partial g_\eps \circ h|/|g_\eps \circ h|,
\]
or
\[
2|\alpha \circ h|^2  \ge |\partial\alpha \circ h|  + |\alpha \circ h|  |\partial g_\eps \circ h|/|g_\eps \circ h|.
\]

The set $\Lambda$ is open, and thus $\partial\Lambda \cap \Lambda = \emptyset$ and also $h(\partial\Lambda) \cap \Omega = \emptyset$.
Moreover, since $h(\partial\Lambda) \subset h(\overline\Lambda) \subset \ol{h(\Lambda)}$, we get
$h(\partial\Lambda) \subset \ol\Omega \setminus \Omega = \pa \Omega$.

Since $w$ is bounded on $\Omega$, the function $v=w \circ h$ is bounded on $\Lambda$, and using the calculations above and
Theorem \ref{thm:halfplane} with condition  (\ref{eq:beta}), we see that 
\[
\sup_{z \in \Lambda} |v(z)| = \sup_{z \in \pa\Lambda} |v(z)|.
\]
Since $h(\Lambda)=\Omega$, $\sup_{z \in \Lambda} |v(z)|=\sup_{z \in \Omega}|w(z)|$. Moreover, since 
$h(\pa\Lambda) \subset \pa\Omega$, we have also
\[
\sup_{z \in \pa\Lambda} |v(z)| \le \sup_{z \in \pa\Omega} |w(z)|.
\]
It follows that $\sup_{z \in \Omega}|w(z)| \le \sup_{z \in \pa\Omega}|w(z)|$ and we obtain equality.

\endpf





We now provide a generalization of the three-lines theorem of Hadamard (see, for example \cite[Thm. 9.4.8]{krantz}
for the classical formulation with $\alpha=0$).

\begin{thm}
Suppose that $a$ and $b$ are real numbers with $0<a<b$, and let $\Omega =\{z \in \CC: a < \re z < b\}$.
Suppose that 
$w$ is a continuous bounded function on $\overline \Omega$
such that (\ref{eq:dbar}) holds in $\Omega$ where $\alpha$ is a $C^1$ function
satisfying 
\beq\label{eq:3linescon}
2|\alpha|^2 \ge |\partial\alpha|+\frac{|\alpha| |\log(M(a)/M(b))|}{b-a}+|\alpha||\partial h_\eps|/|h_\eps|
\eeq
for each $\eps>0$.
Then the function $M$ defined on $[a,b]$ by
\[
M(x)=\sup_{y \in \RR} |w(x+iy)|
\]
satisfies, for all $x \in (a,b)$,
\[
M(x)^{b-a} \le M(a)^{b-x} M(b)^{x-a}.
\]
That is, $\log M$ is convex on $(a,b)$.
\end{thm}

\beginpf
Consider the function $g$ defined on $\overline\Omega$ by
\[
h(z)=M(a)^{(z-b)/(b-a)}M(b)^{(a-z)/(b-a)},
\]
where quantities of the form $M^\omega$ are defined for $M>0$ and $\omega \in \CC$ as $\exp(\omega \log M)$,
taking the principle value of the logarithm.

Now $v:= hw$ satisfies $|v(z)| \le 1$ for $z \in \partial\Omega$, since $|h(a+iy)|=1/M(a)$ and
$|h(b+iy)|=1/M(b)$.

Given that $\overline\partial w=\alpha \overline w$ and that $h$ is holomorphic, then, as we have seen,
$\overline\partial(hw)=\beta \overline{hw}$, where $\beta=\alpha h/\overline h$. Moreover,
$\partial \beta = \partial(\alpha h)/\overline h= (\partial \alpha)(h/\overline h) + \alpha (\partial h)/\overline h$.

Now $\log h= \frac{z-b}{b-a} \log M(a) + \frac{a-z}{b-a}\log M(b)$, and so
\[
\left| \frac{\partial h}{h} \right| = \frac{|\log M(a)/M(b)|}{b-a}.
\]
Thus the condition (\ref{eq:3linescon}) on $\alpha$ implies that
$\beta$ satisfies $2|\beta|^2 \ge | \partial \beta|+
|\beta| |\partial h_\eps|/|h_\eps|$. Hence we can apply Theorem \ref{thm:halfplane} to $v$, and the result follow.

\endpf

\begin{rem}
{\rm  As in Example~\ref{ex:2.2},  rescaling $z$ is helpful here, since if $z$ is reparametrized as $\lambda z$, then
$\partial \alpha$ is divided by $\lambda$ and $b-a$ is also divided by $\lambda$: thus the inequality
(\ref{eq:3linescon})
becomes easier to satisfy.
}
\end{rem}

\section{Weights depending on one variable}

We look at two cases here, for functions defined on a subdomain of $\CC_+$, namely weights $\alpha=\alpha(x)$ and
radial weights $\alpha=\alpha(r)$. We revisit Theorem~\ref{thm:halfplane}.

Since we now have $\partial \alpha=\alpha'/2$, we obtain the following corollary.

\begin{cor}\label{cor:halfplane-x}
Suppose that $\Omega \subset \CC_+$ (not necessarily bounded) and that 
$w$ is a continuous bounded function on $\overline \Omega$
such that (\ref{eq:dbar}) holds in $\Omega$ where $\alpha=\alpha(x)$ is a $C^1$ function satisfying $2|\alpha|^2 \ge | \alpha'|/2+
|\alpha| |\partial h_\eps|/|h_\eps|$ for all $\eps>0$.
Then 
$|w(z)| \le \sup_{\zeta \in \partial \Omega}|w(\zeta)|$ for all $z \in \Omega$.
\end{cor}

Likewise, in polar coordinates $(r,\theta)$ we have 
\[
\partial = \frac{1}{2} \left( e^{-i\theta}\partial_r-\frac{i e^{-i\theta}}{r}\partial_\theta\right),
\]
giving the following result.

\begin{cor}\label{cor:halfplane-x}
Suppose that $\Omega \subset \CC_+$ (not necessarily bounded) and that 
$w$ is a continuous bounded function on $\overline \Omega$
such that (\ref{eq:dbar}) holds in $\Omega$ where $\alpha=\alpha(r)$ is a $C^1$ function satisfying $2|\alpha|^2 \ge | \alpha'|/2+
|\alpha| |\partial h_\eps|/|h_\eps|$ for all $\eps>0$.
Then 
$|w(z)| \le \sup_{\zeta \in \partial \Omega}|w(\zeta)|$ for all $z \in \Omega$.
\end{cor}

Suppose now that $\alpha(x)=ax^\mu$. The condition we require is then
\[
2|a|^2x^{2\mu} \ge |a\mu| x^{\mu-1}/2 + |a|x^\mu \frac{\eps}{|1+\eps z|},
\]
which is only possible for $\mu=-1$. However, it is easy to write down polynomials in $x$ that do not
vanish at $0$ but which satisfy the conditions of Corollary~\ref{cor:halfplane-x}.

\subsection*{Acknowledgments.} The authors are grateful to Joseph Burrier for his assistance. They also thank the referee for some
useful comments.


\begin{thebibliography}{10}

\bibitem{BLRR}
L. Baratchart, J. Leblond, S. Rigat and E. Russ, 
Hardy spaces of the conjugate Beltrami equation. 
{\em J. Funct. Anal.} 259 (2010), no. 2, 384--427.

\bibitem{bers}
L. Bers,  
{\em Theory of pseudo-analytic functions}. Institute for Mathematics and Mechanics, New York University, New York, 1953.

\bibitem{carl}
S. Carl,  
A maximum principle for a class of generalized analytic functions.
{\em Complex Variables Theory Appl.} 10 (1988), no. 2--3, 153--159. 

\bibitem{chen}
 S.-S. Chen,   On a class of quasiconformal functions in Banach spaces. {\em Proc. Amer. Math. Soc.} 37 (1973), 545--548. 

\bibitem{fl}
Y. Fischer and J. Leblond, Solutions to conjugate Beltrami equations and approximation in generalized Hardy spaces. {\em Adv. Pure Appl. Math.} 2 (2011), no. 1, 47--63. 

\bibitem{flps}
Y. Fischer, J. Leblond, J.R. Partington and E. Sincich, Bounded extremal problems in Hardy spaces for the conjugate Beltrami equation in simply-connected domains. {\em Appl. Comput. Harmon. Anal.} 31 (2011), no. 2, 264--285.

\bibitem{IM01}
T. Iwaniec and G. Martin, {\em Geometric function theory and non-linear analysis.} Oxford Mathematical Monographs. The Clarendon Press, Oxford University Press, New York, 2001.

\bibitem{krantz}
S.G. Krantz, { \em Geometric function theory. Explorations in complex analysis.} Cornerstones. Birkh\"auser Boston, Inc., Boston, MA, 2006. 

\bibitem{kravchenko}
V.V. Kravchenko, {\em Applied pseudoanalytic function theory}. With a foreword by Wolfgang Sproessig. Frontiers in Mathematics. Birkh\"auser Verlag, Basel, 2009. 

\bibitem{titchmarsh}
E.C. Titchmarsh, {\em The theory of functions}. Oxford University Press, London, 2nd edition, 1939.

\bibitem{vekua}
I.N. Vekua, {\em Generalized analytic functions}. Pergamon Press, London--Paris--Frankfurt; Addison-Wesley Publishing Co., Inc., Reading, Mass. 1962. 

\end{thebibliography}
\end{document}